\documentclass{article}

\usepackage{amssymb}
\usepackage{graphicx}
\usepackage{lineno}

\newtheorem{theorem}{Theorem}[section]
\newtheorem{lemma}[theorem]{Lemma}

\newtheorem{conjecture}[theorem]{Conjecture}
\newtheorem{problem}[theorem]{Problem}
\newcommand{\DONOTTEX}[1]{}

\newcommand{\notinblk}[1]{}

\makeatletter
\renewcommand{\section}{\@startsection{section}{1}{0em}{2\baselineskip}{1\baselineskip}{\centering\bfseries\Large}}
\renewcommand{\@seccntformat}[1]{\csname the#1\endcsname . \hspace{0.2em}}
\renewcommand{\l@section}{\@dottedtocline{1}{0em}{2.3em}}
\renewcommand{\@evenfoot}{\hfil \textrm{\thepage} \hfil \makebox[0cm][r]{\small\textrm{\today}}}
\renewcommand{\@oddfoot}{\@evenfoot}
\makeatother

\begin{document}


\begin{center}
  {\bf {\LARGE Minimal Connectivity}}\\[7mm]
  {\sc Matthias Kriesell}
\end{center}

\tableofcontents

\bigskip

\begin{center}
  \parbox{10cm}{\it
  A $k$-connected graph such that deleting any edge / deleting any vertex / contracting any edge
  results in a graph which is not $k$-connected is called {\em minimally / critically / contraction-critically $k$-connected}.
  These three classes play a prominent role in graph connectivity theory, and we give a brief introduction
  with a light emphasis on reduction- and construction theorems for classes of $k$-connected graphs.}\\[10mm]
\end{center}

\section{Introduction}
\label{Sintro}

One of the main concerns of graph connectivity theory is to
find reduction methods for classes of $k$-connected graphs.
Such methods can be used for induction proofs, but have also a constructive
counterpart which might be helpful for generating the respective classes.
They also can be employed for setting up problem solving strategies for general graphs:
Roughly, if a graph has a small cutset then we can split it into two smaller parts,
solve the problem, and combine the solutions, and if not then we might reduce 
and use structural properties of $k$-connected graphs.%
\footnote{An instructive example of this is Carsten Thomassen's
brilliant proof for Kuratowski's Theorem \cite{Thomassen1981b}.}

Tutte was the first who studied such methods systematically for $2$- and
$3$-connected graphs. Let us start, as an initial example, with the following version
of his celebrated {\em Wheel Theorem} \cite{Tutte1961}:

\begin{theorem}
  \notinblk{\cite{Tutte1961}}
  \label{Twheel}
  Every $3$-connected simple graph non-isomorphic to a wheel can be
  reduced to a smaller $3$-connected simple graph by either deleting or contracting an edge.%
  \footnote{All graphs considered here are, for the time being, finite and undirected.
  For terminology not defined here, I refer to \cite{BondyMurty2008} and \cite{Diestel2005}.
  To {\em contract} an edge means to delete it and to identify its endpoints; in general,
  that may cause multiple edges.
  A {\em wheel} is obtained from a cycle of length at least $3$
  by adding a new vertex and connecting it to all others by a single edge.}
\end{theorem}

This is a reduction theorem for the class of simple $3$-connected graphs.
It tells us that unless such a graph belongs to a simple subclass of
{\em basic graphs} (the wheels),
we can reduce it to a smaller one by performing a
{\em short  sequence of elementary reductions} (deleting or contracting an edge).
Throughout Section \ref{Sredthms},
we will discuss a number of similar methods for more classes of $k$-connected graphs.
Not all of them will operate in terms of deletion and contraction,
but there is a good reason to concentrate on these:
Contracting or deleting edges, as well as deleting vertices,
will not create minors in the reduced graph which have not been there before,
that is, we will stay inside any given class described by forbidden minors;
the only reason for which the reduction might fail is that the resulting graph is not $k$-connected.

The methods to actually find reducible objects
(for example a single contractible or deletable edge, as in the Wheel Theorem)
are of course encoded in the respective proofs.
They often have a potential for generalization, which would typically answer questions like:
(a) Are there many reducible objects?
(b) Are they all over the place?
(c) Can we reduce in such a way that some additional property survives?
In very basic terms, these are questions on the number and distribution of reducible objects.
In (c), an {\em additional property} might be {\em simplicity} as in the Wheel Theorem;
another such property might be to {\em contain a specific graph as a minor},
and this could lead to what is commonly called a {\em splitter theorem}.
As an example, let us mention Negami's splitter theorem for $3$-connected graphs \cite{Negami1982}:

\begin{theorem}
  \notinblk{\cite{Negami1982}}
  \label{Tsplitter}
  Let $H$ be any $3$-connected graph non-isomorphic to a wheel.
  Then every $3$-connected simple graph non-isomorphic to $H$ and
  containing $H$ as a minor can be reduced to a smaller $3$-connected simple graph
  containing $H$ as a minor by deleting or contracting a single edge.
\end{theorem}  

The results of Section \ref{Sredthms} can be considered as high-end outcomes of a flourishing study
of the distribution of reducible objects in a graph. We start, in Section \ref{Smin}, 
with a brief treatise of {\em minimally $k$-connected graphs},
i. e. $k$-connected graphs where {\em deleting any edge} produces a graph which
is no longer $k$-connected. 
These classes are far too big to become reasonably primitive base classes for a reduction theorem,
even in the case that $k=2$. 
The same applies to {\em critically $k$-connected graphs} studied in Section \ref{Scrit},
i. e. $k$-connected graphs where {\em deleting any vertex} produces a graph
which is no longer $k$-connected,
and  even to  the intersection class of both minimally and critically $k$-connected graphs%
\footnote{For example, every $k$-connected graph where every every edge is incident
with a vertex of degree $k$ and every vertex is adjacent to a vertex of degree $k$ ---
in particular, $k$-regular graphs ---  are among these.}.
Therefore, it is necessary to look at yet another elementary operation --- and see
how far that would lead. According to our above remarks, edge contraction is a good choice.
The corresponding class of {\em contraction-critically $k$-connected graphs},
i. e. $k$-connected graphs where contracting any edge produces a graph which
is no longer $k$-connected, is considered in Section \ref{Sconcrit}.
It is a proper subclass of the critically $k$-connected graphs (however, both criticity concepts can
be treated as special cases of a more general one recalled in Section \ref{Sgeneral}).
For $k \leq 3$, it consists just consists of $K_{k+1}$, which is certainly a sufficiently basic graph.
For $k \geq 4$, there are infinitely many contraction-critically $k$-connected graphs,
and they are rich in a sense. Still, for $k=4$, we may either describe most of them
as line graphs of a class of cubic graphs which in turn admits a constructive characterization.
Alternatively, these contraction critically $4$-connected graphs can be reduced
by contacting {\em two} edges to obtain a smaller $4$-connected graph \cite{Kriesell2002}.
This suggests an idea of what could be true for $5$-connected graphs \cite{Kriesell2002}:

\begin{conjecture}
  \notinblk{\cite{Kriesell2002}}
  \label{C5con}
  There exist integers $b,h$ such that every $5$-connected graph on at least $b$ vertices
  can be reduced to a smaller one by contracting at most $h$ edges.  
\end{conjecture}

As we have pointed out, the corresponding result for $4$-, $3$-, $2$-, and $1$-connected graphs is true.
The reason for Conjecture \ref{C5con} being an interesting question is, however, not, that it would settle
the {\em next} open case, but that it would settle the {\em last} open case: The corresponding
statements for $6$-, $7$-, $\cdots$ connected graphs {\em are not true}, i.e. for $k \geq 6$,
there are minimally $k$-connected graphs showing that we have to contract arbitrarily many edges
to obtain a smaller $k$-connected graph (equivalently, the gap between the order of such a graph
and the maximum order of its lower neighbors in the minor relation, restricted to $k$-connected graphs,
is arbitrarily large) \cite{Kriesell2002}.

This indicates that it is maybe worthwhile to look at larger substructures of the graph than just
vertices or edges and to study possible ways to employ them for reduction. This is the topic of
Sections \ref{Shcks} and \ref{Spart}.

Many problems involving vertex-connectivity have a literal counterpart in terms of
edge-connectivity. Although in most cases these turn out to be easy if not trivial,
the restriction of the original vertex-connectivity problems to {\em line graphs}
is a little bit more demanding and may serve as a touchstone, as it is finally illustrated in Section \ref{Slgr}.

\section{Edge Deletion}
\label{Smin}

We start with a brief section on graphs where deleting any edge decreases its connectivity.
Let $G$ be a $k$-connected graph%
\footnote{That is, $G$ has more than $k$ vertices and cannot be separated by removing less than $k$ vertices,
i.e. $G-X$ is connected for all $X \subseteq V(G)$ with $|X|<k$. The smallest $k$ for which $G$ is
$k$-connected is the {\em connectivity} of $G$, denoted by $\kappa(G)$.}.
An edge $e$ of $G$ is called {\em $k$-essential} if $G-e$ is not $k$-connected,
and $G$ is called {\em minimally $k$-connected} if every edge is $k$-essential.
The connectivity of such a graph is equal to $k$.
Minimally $k$-edge-connected graphs are very well understood which is mostly
due to the following fundamental statement relating vertices of
degree $k$ and $k$-essential edges \cite[Satz 1]{Mader1972}.

\begin{theorem}
  \notinblk{\cite[Satz 1]{Mader1972}}
  \label{Tess}
  Let $G$ be a $k$-connected graph and $C$ be a cycle of $k$-essential edges in $G$.
  Then there exists a vertex in $V(C)$ of degree $k$ in $G$.
\end{theorem}

As an immediate consequence, in a minimally $k$-connected graph $G$ the vertices
of degree larger than $k$ induce a forest \cite{Mader1972},
and it is easy to derive that such a graph must have more than $\frac{k-1}{2k-1}|V(G)|$
many vertices of degree $k$ \cite{Mader1972}.
Moreover, Theorem \ref{Tess} implies that every 
minimally $k$-connected graph $G$ is the edge disjoint union of a tree and $k-1$ forests%
\footnote{The proof by induction on $k$ starts obviously for $k=1$, whereas for $k>1$
$G$ contains a minimally $(k-1)$-connected subgraph $H$, for which the statement is true; if there was
a cycle consisting of edges in $E(G)-E(H)$ then, by Theorem \ref{Tess}, it contains a vertex $x$
of degree $k$ in $G$ --- but then $x$ had degree $k-2$ in $H$, a contradiction.}.
In particular, the average degree of a minimally $k$-connected graph is less than $2k$.

Another consequence of Theorem \ref{Tess} is that every minimally $k$-connected graph has
at least $k+1$ vertices of degree $k$ \cite{Mader1972}.
This implies that --- unlike in the case of vertex deletion (see the next section) ---
deleting an edge $e$ from a $k$-connected graph $G$, where $k \geq 2$,
will almost never produce a minimally $(k-1)$-connected graph:
The only exception is that $G$ is a cycle.%
\footnote{If $G-e$ is minimally $(k-1)$-connected then it has $k$ vertices of degree $k-1$.
But on the other hand it has at most two vertices of degree $k-1$, namely the endvertices of $e$.
Hence $k=2$ and $G-e$ is minimally $1$-connected with two vertices of degree $1$,
so that $G$ must be a cycle.}

Halin proved earlier that every triangle of $k$-essential edges in a $k$-connected graph
contains at least {\em two} vertices of degree $k$ (see \cite[Satz 3']{Halin1969}).

For a comprehensive survey on minimally $k$-connected graphs and digraphs
see \cite{Mader1993}.

\section{Vertex Deletion}
\label{Scrit}


Let $G$ be a $k$-connected graph.
A vertex $x$ of $G$ is called {\em $k$-essential} if $G-x$ is not $k$-connected,
and $G$ is called {\em critically $k$-connected} if every vertex is $k$-essential.
Again, in this case, $k$ must be equal to $\kappa(G)$.
Any $k$-connected graph where every vertex is adjacent to a vertex of degree $k$ is critically $k$-connected.
Clearly, such graphs might have vertices of arbitrarily large
degree unless $k=1$. One might ask if, as for minimally $k$-connected graphs,
a critically $k$-connected graph must have vertices of degree $k$.
This turns out to be true for $k=2$ and $k=3$, by the following Theorem ---
but it is wrong in general. To see this, observe that if $G$ is a critically $k$-connected graph then
(a) the lexicographic product%
\footnote{The {\em lexicographic product} $G[H]$ of two simple loopless graphs $G,H$
is the graph defined by $V(G[H]):=V(G) \times V(H)$ and
$E(G[H]):=\{(x,a)(y,b):\,xy \in E(G) \vee (x=y \wedge ab \in E(H))\}$.
Roughly, the vertices of $G$ are replaced by disjoint copies of $H$,
and two such copies are connected by all possible edges if their corresponding vertices
in $G$ are adjacent.}
$G[K_\ell]$ of $G$ and a complete graph of order $\ell$
is critically $k \ell$-connected with $\delta(G[K_\ell])=\delta(G) \cdot \ell + \ell-1$, and
(b) the graph $G*K_1$ obtained from $G$ by adding a new vertex and making it adjacent
to all others is critically $(k+1)$-connected with $\delta(G*K_1)=\delta(G)+1$.
Hence given $k \geq 2$ and any critically $2$-connected graph $G$, the graphs
$G[K_{k/2}]$ for even $k$ and $G[K_{(k-1)/2}]*K_1$ for odd $k$ are critically $k$-connected
of minimum degree $\lfloor \frac{3}{2} k -1 \rfloor$. At the same time, they are sharpness
examples of the following theorem from \cite{Lick1968}.

\begin{theorem}
  \notinblk{\cite{Lick1968}}
  \label{T1}
  For every critically $k$-connected graph $G$,
  \[ \textstyle \delta(G) \leq \frac{3}{2} \cdot k -1. \]
\end{theorem}

This is a consequence of a more general result.
To describe it, we need to introduce some terminology.
For any graph $G$, let us denote by ${\cal T}(G)$ the set of all separating
vertex sets of cardinality $\kappa(G)$. A {\em fragment} of $G$ is
the union of the vertex sets of at least one but not all components
of $G-T$, for some $T \in {\cal T}(G)$. For a non-complete graph,
fragments of minimum cardinality
are called {\em atoms}, and this cardinality is denoted by $a(G)$.
It is easy to see that a vertex of a non-complete graph $G$ of connectivity $k$ is $k$-essential
if and only if it is contained in some member of ${\cal T}(G)$; consequently, $G$
is critically $k$-connected if and only if $\bigcup {\cal T}(G)=V(G)$.
If $A$ is an atom then the neighbors of any $x \in A$
are in $(A-\{x\}) \cup N_G(A)$ and $|N_G(A)|=\kappa(G)$ because $N_G(A) \in {\cal T}(G)$,
and so 
\[ \textstyle \delta(G) \leq \kappa(G)+a(G)-1. \]
Hence Theorem \ref{T1} is an immediate consequence of the following from the following Theorem in \cite{Mader1971}:

\begin{theorem}
  \notinblk{\cite{Mader1971}}
  \label{T2}
  For every non-complete critically $k$-connected graph $G$,
  \[ \textstyle a(G) \leq k/2. \]
\end{theorem}

This is an almost immediate consequence of the following fundamental lemma due to Mader \cite{Mader1971}.

\begin{lemma}
  \notinblk{\cite[Satz 1]{Mader1971}}
  \label{L1}
  If $A$ is an atom of a graph $G$ and $T \in {\cal T}(G)$ contains at least one vertex from $A$ then
  (i) $A \subseteq T$ and (ii) $|A| \leq \frac{|T-N_G(A)|}{2}$.%
\footnote{This is a good justification for the term {\em atom},
from greek $\acute{\alpha}\tau{o}\mu{o}\varsigma$, meaning uncuttable/indivisible,
which is resembled in in (i) --- plus the connotation of being small, resembled in (ii).}%
\end{lemma}

Later, we will come to a yet more general version of this (Lemma \ref{L1g}).

Unlike for minimally $k$-connected graphs, the average degree of critically $k$-connected graphs
cannot be bounded by a function of $k$ \cite{Kriesell2006b}.

Lemma 1 implies that distinct atoms of a $k$-connected graph are disjoint.
This was first proved by Watkins, and he used
it to relate degree and connectivity of vertex transitive graphs
(see \cite{Watkins1970}, where there is also an independent proof of Theorem \ref{T1},
for vertex transitive graphs). The question of relating degree and connectivity of a vertex transitive graph
goes back to Vizing (see page 130 of \cite{Vizing1968}, and be aware that {\em vertex transitive} graphs are called {\em regular} there).
For even $k$, the lexicographic products $C_\ell[K_{k/2}]$ show that the bound
of Theorem \ref{T1} is sharp for vertex transitive graphs,
but for odd $k$ our examples $G[K_{(k-1)/2}] * K_1$ are not vertex transitive.
Indeed, in this case we can improve the bound as follows.
First, (i) of Lemma \ref{L1} implies that distinct atoms of a critically $k$-connected graph are disjoint.
Moreover, if $G$ is vertex transitive then the atoms form a system of imprimitivity and,
again by (i) of Lemma \ref{L1}, every $T \in {\cal T}(G)$ is the disjoint union of atoms.
Therefore, $a(G)$ divides $|T|=\kappa(G)$, and $a(G) \not= |T|$
by (ii) of Lemma \ref{L1}. Hence $a(G) \leq k/p$, where $p$ is the smallest prime divisor of $k$,
which improves the bound from Theorem \ref{T1} to
\[ \textstyle \delta(G) \leq \frac{p+1}{p} \cdot k - 1. \]
This is sharp because of the vertex transitive graphs $G[K_{k/p}]$, where $G$ is
any vertex transitive $p$-regular graph of connectivity $p$.
See also Section 4 in \cite{Watkins1970}.
Later, Jung \cite{Jung1976} used similar methods to analyze graphs where every vertex is contained
in the same number of {\em $k$-atoms}, i. e. fragments $F$ with $|F| \geq k$ and
$|V(G)-(F \cup N_G(F))| \geq k$ of minimum cardinality. In particular,
this lead to an interesting structural classification of vertex transitive graphs in terms
of fragment clusters \cite{Jung1976}.

It has been observed by Watkins \cite{Watkins1970} that
the disjointness of distinct atoms implies $a(G)=1$ (and hence
$\delta(G)=\kappa(G)$) for an {\em edge-transitive graph} $G$.
Here is his argument: Assume, to the contrary, that $a(G) > 1$, and let $A$ be an atom.
Then there exists a path $xyz$ with $x \in N_G(A)$ and $y,z \in A$,
and an automorphism of $G$ which maps $yz$ to $xy$.
Clearly, $xy$ is contained in an atom distinct but not disjoint from $A$, contradiction.

Maurer and Slater \cite{MaurerSlater1977}
suggested to generalise the concept of critically $k$-con\-nec\-ted graphs.
They defined a graph $G$ to be {\em $\ell$-critically $k$-connected} if $\kappa(G-X)=k-\ell$
for all $X \subseteq V(G)$ with $|X| \leq \ell$; equivalently, $\kappa(G)=k$ and either $G \cong K_{k+1}$ or every
$X \subseteq V(G)$ with $|X| \leq \ell$ is a subset of some $T \in {\cal T}(G)$.
Any such $G$ is called {\em $\ell$-critical}.
Obviously, $G$ is critically $k$-connected if and only if $G$ is $1$-critically $k$-connected,
and every $\ell$-critically $k$-connected graph is $\ell'$-critically $k$-connected whenever $1 \leq \ell' \leq \ell$.
Trivially, $\ell \leq k$, but Maurer and Slater conjectured that even $\ell \leq k / 2$ holds
unless $G$ is isomorphic to $K_{k+1}$. This has first been proved by Su
(see \cite{Su1988} from 1988, and, easier to access, \cite{Su1993} from 1993):

\begin{theorem}
  \notinblk{\cite{Su1993}}
  \label{T3}
  Suppose there exists an $\ell$-critically $k$-connected graph non-iso\-mor\-phic to $K_{k+1}$.
  Then $\ell \leq k / 2$.
\end{theorem}

The graphs 
\[ \textstyle S_\ell:=K_{\ell+1}[\overline{K_2}] \]
and $S_\ell-x$ are $\ell$-critically $2\ell$-connected and $(\ell-1)$-critically $(2\ell-1)$-connected, respectively,
thus showing that the bound in Theorem \ref{T3} is sharp.

Later, in 1998, Jord\'an found a very elegant argument of Theorem \ref{T3},
which is one of the pearls of this part of the theory \cite{Jordan1998}.
It depends on Theorem \ref{Tess}.
So let $G$ be an $\ell$-critically $k$-connected graph, and let $A$ be an atom of $G$.
By Lemma \ref{L1} it follows easily that $G-A$ is $(\ell-1)$-critically $(k-|A|)$-connected,
and that every fragment $F$ of $G-A$ is a fragment of $G$, where $N_G(F)=N_{G-A}(F) \cup A$.
Each fragment of $G-A$ contains a vertex of $S:=N_G(A)$. Therefore, we may add
some set $N$ of new edges between vertices of $S$ such that $(G-A)+N$ is $(k-|A|+1)$-connected.
If we take $N$ inclusion minimal with respect to this property then every new edge 
will be $(k-|A|+1)$-essential for $(G-A)+N$, so that $N$ forms a forest on $S$ by
Theorem \ref{Tess}.
The edges of a forest (of {\em any} bipartite graph) can be covered by at most half of its vertices,
so that there exists a set $X$ of at most $|S|/2$ vertices in $S$ meeting everybody from $N$.
As for every smallest separating set of $G-A$, there must be an edge from $N$
connecting two distinct components of $G-A$ (for otherwise $\kappa((G-A)+N)=k-|A|$),
we deduce that $X$ is not contained in a smallest separating set of $G-A$.
Consequently, $|X| \geq \ell$, as $G-A$ is $(\ell-1)$-critically $(k-|A|)$-connected.
Therefore, $\ell \leq |X| \leq |S|/2=\kappa(G)/2$, as desired.

Mader conjectured several properties of $\ell$-critically $k$-connected graphs stronger
than the statement of Theorem \ref{T3} in a survey paper from 1984 \cite{Mader1984},
which are, meanwhile, all proved.
The possibly most difficult result along these lines is the following theorem from \cite{Su2004}:

\begin{theorem}
  \notinblk{\cite{Su2004}}
  \label{T4}
  If $G \not\cong K_{k+1}$ is $\ell$-critically $k$-connected then $G$ contains $2\ell+2$ disjoint fragments.
\end{theorem}

This has been proved by Su, too \cite{Su2004}, and his proof is very ingenious and difficult.
Unlike in the case of Theorem \ref{T3}, no simpler proof has been found so far.
(Refining his method outlined above, Jord\'an gave a --- still simple --- proof
that there is an antichain\footnote{With respect to $\subseteq$.} of $2\ell+2$ fragments \cite{Jordan1998}.) ---
Given an atom $A$ of an $\ell$-critically $k$-connected graph $G$, we may
apply Theorem \ref{T4} to the $(\ell-1)$-critically $(k-|A|)$-connected graph $G-A$ (see above).
It is easy to see that every fragment of $G-A$ contains at least $|A|$ vertices from $S:=N_G(A)$,
so that $(2 (\ell-1)+2) \cdot |A| \leq |S|=k$, which yields the following generalization of
Theorem \ref{T2} (also conjectured by Mader \cite{Mader1984}); see \cite{Su2004}:

\begin{theorem}
  \notinblk{\cite{Su2004}}
  \label{T4a}
  If $G \not\cong K_{k+1}$ is $\ell$-critically $k$-connected then $a(G) \leq \frac{k}{2 \ell}$.
\end{theorem}

The graphs $S_\ell[K_m]$ are $\ell$-critically $2 \ell m-connected$, and they
have exactly $2\ell+2$ fragments, which are pairwise disjoint and have cardinality $m$ ---
hence the bounds in Theorem \ref{T4a} (and Theorem \ref{T4}) are sharp.

Skew to these Theorems, there is the following statement on the
{\em extremely critically connected graphs},
that is, the $\ell$-critically $2\ell$-connected graphs (another former conjecture by Mader \cite{Mader1984}).
See \cite{Kriesell2000}, \cite{Kriesell2006}, and \cite{SuYuanZhao2003}:

\begin{theorem}
  \notinblk{\cite{Kriesell2000} \cite{Kriesell2006} \cite{SuYuanZhao2003}}
  If $G \not\cong K_{2\ell+1}$ is $\ell$-critically $2\ell$-connected and $\ell \geq 3$ then
  $G \cong K_{\ell+1}[\overline{K_2}]=S_\ell$.
\end{theorem}

A substantial part of the proof is considered with the small cases $\ell \in \{3,4\}$ \cite{Kriesell2000},
and, in \cite{Kriesell2006}, to the cases up to approximately $\ell=20$. For larger values, the
statement turned out to be {\em less and less difficult} in a sense \cite{Kriesell2006}.
The proof in \cite{Kriesell2006} has been
developed while the statement of Theorem \ref{T4} was still an open question, and is, thus,
independent from Theorem \ref{T4}. However, it is possible to simplify part of the
work for the small cases by using Theorem \ref{T4}, as it has been demonstrated in \cite{SuYuanZhao2003}.
A characterization of the extremely critically connected graphs for {\em odd} connectivity $\geq 5$,
that is, the $\ell$-critically $(2\ell+1)$-connected graphs, might be achieved
in the future.

There are several open questions on $\ell$-critically $k$-connected graphs,
also for small values of $\ell$. For the first one, see \cite{Mader1984}:

\begin{conjecture}
  \notinblk{(See \cite{Mader1984})}
  \label{C1}
  Is there an $\ell$ such that every $\ell$-critically $k$-connected graph contains a $K_4$?
\end{conjecture}

The octahedron $S_2=K_3[\overline{K_2}]=K_{2,2,2}$ is $2$-critically $4$-connected and does not contain $K_4$,
thus showing $\ell \geq 3$ if $\ell$ as in Conjecture \ref{C1} existed. 

If $a,b$ are vertices of a $3$-critically $k$-connected graph and $A$ is an atom
of the $1$-critically $(k-2)$-connected graph $G-\{x,y\}$, then there exist neighbors
$c,d \in A$; by Theorem \ref{T2}, $|A| \leq (k-2)/2$, which implies that $c,d$ have
a common neighbor in $A \cup N_G(A)$. This shows that $a,b$ are at distance at most $4$,
and as they have been chosen arbitrarily, $G$ has diameter at most $4$ \cite{Mader1977}.
Mader asked if there exists any $3$-critically $k$-connected graph of diameter $4$ or $3$ \cite{Mader1984}.
In \cite{Kriesell2007},
I constructed for every $\ell \geq 3$ an $\ell$-critically $k$-connected graphs of diameter $3$
(the smallest one is $122$-connected and  has $252$ vertices).
I doubt that there are any of diameter $4$ (cf. \cite{Kriesell2007}):

\begin{conjecture}
  \notinblk{(See \cite{Kriesell2007}.)}
  Every $3$-critically $k$-connected graph has diameter at most $3$.
\end{conjecture}

In order to extend the construction from \cite{Kriesell2007} to produce $\ell$-critically $k$-connected
graphs of diameter $4$, it would be necessary to find sufficiently large $\ell$-critically
$k'$-connected graphs. This turned out to be impossible for $\ell \geq 5$ \cite{Kriesell2007}, and
it would be impossible also for $\ell \in \{3,4\}$ if the answer to the following question
\cite[Conjecture 4.3]{Mader2002} (see also \cite{Mader1977}) was affirmative:

\begin{conjecture}
  \notinblk{\cite[Conjecture 4.3]{Mader2002} (see also \cite{Mader1977})}
  There exists a $c>0$ such that every $3$-critically $k$-connected graph has at most
  $ck^{3/2}$ many vertices.
\end{conjecture}

In \cite[Corollary 1.3]{Mader2002} it has been shown that every
$3$-critically $k$-connected graph has at most
$2k^2 - k$ many vertices, improving the constant in the bound $6k^2$ from \cite[Satz 2]{Mader1977}.
It is not known whether there exists a $c>0$ and a $\lambda<2$ such that every $3$-critically
$k$-connected graph has at most $ck^\lambda$ vertices. It is not even known if this
holds for $\ell$-critically $k$-connected graphs for a sufficiently large $\ell$ (not depending on $k$).

Regarding our quest for reduction theorems for a class ${\cal C}$ of $k$-connected graphs,
it is possibly difficult to exploit knowledge on $\ell$-critically $k$-connected graphs for $\ell \geq 3$.
However, a $k$-connected graph $G$ has two vertices $x \not= y$ such that their identification
produces a smaller $k$-connected graph $G'$ unless $G$ is $2$-critically $k$-connected.
These vertices may be adjacent or not; if they {\em are} connected by an edge $e$
then $G'=G/e$, and this will be discussed in Section \ref{Sconcrit}. If
they are not then we would possibly leave class ${\cal C}$ (if it is described by forbidden minors).
However, identifying non-adjacent vertices can be used to reduce $k$-connected
{\em bipartite} graphs (and contraction is useless there!), as it will be explained in Section \ref{Sredthms}.

\section{Edge Contraction}
\label{Sconcrit}

Let again $G$ be a $k$-connected graph.
An edge $e$ of $G$ is called {\em $k$-contractible} if $G/e$ is $k$-connected,
and $G$ is called {\em contraction-critically $k$-connected} if there are no $k$-contractible edges.
Again, in this case, $k$ must be equal to  $\kappa(G)$.
There is a plethora of results on the existence, number, and distribution of $k$-contractible edges,
and the reader is refered to my survey from 2002 \cite{Kriesell2002}. 
Here, I will repeat the basics and concentrate on more recent developments.

Trivially, every edge of a $1$-connected graph non-isomorphic to $K_2$ is $1$-con\-trac\-tible,
and it is easy to see that every vertex in a $2$-connected graph $G \not\cong K_3$
is incident with some $2$-contractible edge. 
Tutte's Theorem \ref{Twheel} from Section \ref{Sintro} implies that every
$3$-connected graph $G \not\cong K_4$ has a $3$-contractible edge (but the wheels of order at least $5$
and many other examples show that there might be vertices not incident with any $3$-contractible edge).
From these statements we deduce:
 
\begin{theorem}
  For $k \leq 3$, $K_{k+1}$ is the only contraction-critically $k$-connected graph.
\end{theorem}

In other words, for $k \leq 3$, every $k$-connected graph
non-isomorphic to $K_{k+1}$ can be reduced to a smaller
$k$-connected graph by contracting {\em one single edge}, 
and $K_{k+1}$ is the only minor-minimal $k$-connected graph.
This has been generalised by Egawa \cite[Theorem B]{Egawa1991} as follows. 

\begin{theorem}
  \notinblk{\cite[Theorem B]{Egawa1991}}
  \label{T4b}
  If $G \not\cong K_{k+1}$ is contraction-critically $k$-con\-nec\-ted then $a(G) \leq k/4$.
\end{theorem}

As every $2$-critically $k$-connected graph is contraction-critically $k$-connected,
Theorem \ref{T4b} also generalises Theorem \ref{T4a}, restricted to the case that $\ell=2$.
In the spirit of Maurer's and Slater's generalization of criticity (see Section \ref{Scrit}),
Mader called a graph $G$ {\em $\ell$-con-critically $k$-connected}
if $\kappa(G-X)=k-\ell$ for all $X \subseteq V(G)$ with $|X| \leq \ell$
and $G[X]$ connected \cite{Mader2002}; equivalently, $\kappa(G)=k$ and either $G \cong K_{k+1}$ or every
$X \subseteq V(G)$ with $|X| \leq \ell$ and $G[X]$ connected is a subset of some $T \in {\cal T}(G)$.
Clearly, every $\ell$-critically $k$-connected graph is
$\ell$-con-critically $k$-connected, but the converse is not true (which is, by
the way, not easy to certify, see \cite[Section 5]{Mader2002}).
Opposed to the situation of $3$-critically $k$-connected graphs (see Section \ref{Scrit}),
it is not possible to bound the number of vertices in a $3$-con-critically  $k$-connected graph
by a function of $k$ \cite{Mader2004}. However, Theorem \ref{T3} generalises as follows \cite{Mader2004}:

\begin{theorem}
  \notinblk{\cite{Mader2004}}
  \label{T3a}
  Suppose there exists an $\ell$-con-critically $k$-connected graph non-isomorphic to $K_{k+1}$.
  Then $\ell \leq k / 2$.
\end{theorem}

There are infinitely many contraction critically $4$-connected graphs.
Fortunately, they are all known, due to the following result by Fontet \cite{Fontet1978,Fontet1979}
and Martinov \cite{Martinov1981,Martinov1982,Martinov1990} (cf. \cite{Mader1984}):

\begin{theorem}
  \label{T5}
  The contraction-critically $4$-connected graphs are the squares of cycles of length at
  least $5$ and the line graphs of cubic cyclically $4$-edge-connected graphs.
\end{theorem}

From this it is easy to see that every $4$-connected graph non-isomorphic to $K_5$ and $K_{2,2,2}$
can be reduced to a smaller $4$-connected graph by contracting {\em one or two edges} (so that
$K_5$ and $K_{2,2,2}$ are the only minor minimal $4$-connected graphs).
As it has been mentioned in the introduction, it is a burning question
if a similar result was true for $5$-connected graphs (see Conjecture \ref{C5con}).

Another, less straightforward way to turn Theorem \ref{T5} into a reduction theorem
is to reduce cubic cyclically $4$-edge-connected graphs, as follows:
To {\em suppress} a vertex $x$ of degree at most $2$ in a graph $G$ means to delete
it and add an edge from $a$ to $b$ if $a,b$ were distinct non-adjacent neighbors of $x$.
To {\em homotopically delete an edge} $e$ in a graph means to delete it and
to suppress its endvertices should they have degree at most $2$. The result
is denoted by $G -- e$. Observe that if $G$ is cubic and triangle-free then $G -- e$
is cubic. Now every simple, cyclically $4$-edge-connected cubic graph 
non-isomorphic to $K_4$ or the skeleton of a $3$-dimensional cube can
be reduced to a smaller cyclically $4$-edge-connected cubic graph by homotopically delete
an edge (see \cite{Mader1984}). In the line graph this means to delete a vertex (of degree $4$) and
contract the two (disjoint) edges in its former neighborhood%
\footnote{The drawback of this, compared to just contracting any triangle in
the line graph, is, that, in general, we cannot do it everywhere.}.
So if $G$ is contraction-critically $4$-connected then we may either reduce it this way,
or $G$ is the square of a cycle of length $\ell \geq 5$. If $\ell=5$ then $G \cong K_5$,
if $\ell=6$ then $G \cong K_{2,2,2}$, and if $\ell \geq 7$ then we may contract two
edges to obtain the square of a cycle of length $\ell-2$.

Like for minimally $k$-connected graphs, and unlike for critically $k$-connected graphs,
the average degree of a contraction-critically $k$-connected graph can be bounded
from above by a function $f$ of $k$ \cite{Kriesell2006b}. In fact, $f(k) \leq c k^2 \log k$
for some constant $c$ \cite{Kriesell2006b}.
As there are contraction-critically $k$-connected graphs of average degree $k^2/6$
\cite{Kriesell2001}, the bound is sharp up to the logarithmic factor --- 
of which I believe that it can be omitted  \cite[Conjecture 2]{Kriesell2006b}:

\begin{conjecture}
  \notinblk{\cite[Conjecture 2]{Kriesell2006b}}
  There exists a constant $c$ such that every finite $k$-connected graph of average degree at least
  $c k^2$ admits a $k$-contractible edge.
\end{conjecture}

In \cite{Kriesell2006b}, I constructed contraction-critically $5$-connected graphs of average degree $10$,
and conjectured that this would be sharp. However, Ando gave an example
of a contraction-critical $5$-connected graph of average degree $12.5$ [personal communication].

A $k$-contractible edge can also be forced by degree sum conditions:
It has been proved in \cite{Kriesell2001} for $k \geq 4$ and $k \not=7$ and
in \cite{SuYuan2005} for $k=7$, that if the sum of the degrees of any pair
of vertices at distance $1$ or $2$ in a $k$-connected graph non-isomorphic to $K_{k+1}$ is at least
$2 \lfloor \frac{5}{4} k \rfloor  -1$ then there exists a $k$-contractible edge.
In fact, even the following weaker degree sum condition turned out to be sufficient
(as I conjectured for all $k$ in \cite{Kriesell2001}), see \cite{SuYuan2006}:

\begin{theorem}
  \notinblk{\cite{SuYuan2006}}
  Let $k \geq 8$.
  If the sum of the degrees of any two adjacent vertices in a $k$-connected graph is at least
  $2 \lfloor \frac{5}{4} k \rfloor  -1$ then there exists a $k$-contractible edge.
\end{theorem}

The degree sum bound is sharp for every $k$ \cite{Kriesell2001}.

\section{Generalised Criticity}
\label{Sgeneral}

Most of the considerations on fragments presented in Section \ref{Scrit}
can be generalised to fragments whose neighborhood contains a member
of a specified set ${\cal S}$ of vertex sets of the graph in question. 
This approach has been worked out by Mader in \cite{Mader1988}.
Given a graph $G$ and ${\cal S} \subseteq {\cal P}(V(G))$,
we call a fragment of $G$ an {\em ${\cal S}$-fragment} if $S \subseteq N_G(F)$ for some
$S \in {\cal S}$. An {\em ${\cal S}$-atom} is a minimum ${\cal S}$-fragment, and its cardinality
is denoted by $a_{\cal S}(G)$. Lemma \ref{L1} generalises as follows \cite{Mader1988}:

\begin{lemma}
  \notinblk{\cite{Mader1988}}
  \label{L1g}
  Let $G$ be a graph and ${\cal S} \subseteq {\cal P}(V(G))$.
  If $A$ is an ${\cal S}$-atom and there exist $S \in {\cal S}$ and $T \in {\cal T}(G)$
  such that $S \subseteq T \cap (A \cup N_G(A))$ and $A \cap T \not= \emptyset$ 
  then $A \subseteq T$ and $|A| \leq |T-N_G(A)|/2$.
\end{lemma}

In all the results on critical graphs above,
the set of critical objects --- let that be vertices, vertex sets, or edges ---
was dense in the rough sense that {\em everywhere} in the graph we could find them.
However, in most basic situations we need only a slightly weaker density condition,
which just ensures that the preconditions of Lemma \ref{L1g} are satisfied:
So let us call a graph $G$ of connectivity $k$ {\em ${\cal S}$-critically $k$-connected},
where ${\cal S} \subseteq {\cal P}(V(G))$,
if ${\cal S} \not= \emptyset$, every $S \in {\cal S}$ is a subset of some $T \in {\cal T}(G)$,
and for every ${\cal S}$-fragment $A$ there exist $S \in {\cal S}$ and $T \in {\cal T}(G)$
such that $S \subseteq T \cap (A \cup N_G(A))$ and $A \cap T \not= \emptyset$.
Using this notion, Lemma \ref{L1g} implies the following theorem from \cite{Mader1988}, by literally the same argument
that led from Lemma \ref{L1} to Theorem \ref{T2}:

\begin{theorem}
  \notinblk{\cite{Mader1988}}
  \label{T2g}
  For every ${\cal S}$-critically $k$-connected graph $G$,
  \[ a_{\cal S}(G) \leq k/2. \]
\end{theorem}
Mader designed this concept as a common generalization of many of the previously mentioned
criticity concepts, and of others. Suppose that $G$ is a non-complete graph. Then:
\begin{enumerate}
  \item
    $G$ is critically $k$-connected iff $G$ is $\{ \{x\}:\,x \in V(G) \}$-critically $k$-connected.
  \item
    $G$ is $\ell$-critically $k$-connected iff $G$ is $\{ X:\,X \subseteq V(G),\,|X| \leq \ell\}$-critically $k$-connected.
  \item
    $G$ is $\ell$-con-critically $k$-connected iff $G$ is $\{ X:\,X \subseteq V(G),\,|X| \leq \ell,\,G[X]$ is connected$\}$-critically $k$-connected.
  \item
    $G$ is contraction-critically $k$-connected iff $G$ is $\{V(e):\, e \in E(G)\}$-critically $k$-connected.
  \item
    $G$ is {\em almost critically $k$-connected} iff $G$ is $\{\emptyset\}$-critically $k$-connected.
  \item
    $G$ is {\em clique-critically $k$-connected} iff
    $G$ is $\{X:\, X \subseteq V(G),\, G[X] \mbox{ complete}\}$-critically $k$-connected.
\end{enumerate}

Let us briefly consider the last two items of the list.
Accordingly, a $k$-connected graph $G \not\cong K_{k+1}$
is almost critically $k$-connected if and only if every fragment contains
a vertex from $\bigcup {\cal T}(G)$. These graphs are important when studying the distribution of contractible edges.
So suppose that a vertex $x$ of some graph $G \not\cong K_{k+1}$ of connectivity $k$
is not incident with a $k$-contractible edge. Then it is easy to see that $G$ is
${\cal S}_x$-critical, where ${\cal S}_x:=\{V(e):\, e \in E_G(x)\}$. Every fragment $F$
of $G-x$ is a fragment of $G$, where $N_G(F)=N_{G-x}(F) \cup \{x\}$ and both $F$ and $N_G(F)$
contain a vertex from $N_G(x)$. 
(In particular, $G-x$ is almost critically $(k-1)$-connected.)
Following Mader's argument from \cite{Mader1988},
let us use this fact to prove that there is a triangle in $G$ `close' to $x$, by considering
an ${\cal S}_x$-atom $A$: If $A$ consists of a single vertex $y$ then $x$ is on a triangle,
formed by $x$, $y$, and any neighbor of $x$ in $N_G(A)$. Otherwise, $A$ must contain a pair
of adjacent vertices, and they must have a common neighbor in $A \cup N_G(A)$ because
$|A| \leq k/2$ by Theorem \ref{T2g} (in fact, $|A| \leq (k-1)/2$ by the appropriate application of Lemma \ref{L1g}).
For a detailed treatise of this interconnection, see \cite{Mader1988}.
Let us just mention that we get one of the main results in \cite{Thomassen1981}, as a corollary \cite{Mader1988}:

\begin{theorem}
  \notinblk{\cite{Thomassen1981} \cite{Mader1988}}
  \label{Ttriangle}
  Every  non-complete triangle free $k$-connected graph contains a $k$-contractible edge.
\end{theorem}

A graph is clique-critical $k$-connected if and only if it is non-complete and $\kappa(G-V(K))=\kappa(G)-|V(K)|$
for every clique $K$ in $G$ (roughly, deleting a clique will always exhaust its potential
of decreasing $\kappa$, that is, decrease it by its order). A tantalizing question is, whether
there are clique-critical graphs at all. Mader conjectured \cite{Mader1988}:

\begin{conjecture}
  \notinblk{\cite{Mader1988}}
  \label{Ccc}
  There is no clique-critical $k$-connected graph.
\end{conjecture}

If this was true, then, for example, every $\ell$-critically $k$-connected graph must contain
a clique on $\ell+1$ vertices, which provides an affirmative answer to Conjecture \ref{Ccc}.
As a clique-critically $k$-connected graph is contraction-critically $k$-connected,
it would contain a triangle by Theorem \ref{Ttriangle} ---
but it is not known whether it would contain a $K_4$, cf. Conjecture \ref{C1}. Mader proved
that Conjecture \ref{Ccc} is true for $k \leq 6$ \cite{Mader1988}.

\section{Reduction Methods}

\label{Sredthms}

Let us summarise some reduction methods for $k$-connected graphs.
The first set runs in terms of $k$-contractible edges; for a comprehensive survey, see \cite{Kriesell2002}.
Trivially, every edge of a $1$-connected graph non-isomorphic to $K_2$ is $1$-contractible.
The respective statement for $2$-connected graphs is not true, as it might happen that the
endvertices of some edge separate. However, we have%
\footnote{As an immediate consequence of Lemma \ref{L1g} applied to ${\cal S}_x$ as in the preceeding paragraph.}:

\begin{theorem}
  \label{Tred2}
  Every vertex of a $2$-connected graph non-isomorphic to $K_3$ is incident with a $2$-contractible edge.
\end{theorem}

For $3$-connected graphs, we have the following, as an immediate consequence either of the
Wheel Theorem (Theorem \ref{Twheel}), or of Theorem \ref{T4b}:

\begin{theorem}
  \notinblk{\cite{Tutte1961}}
  \label{Tred3}
  Every $3$-connected graph non-isomorphic to $K_4$ has a $3$-con\-trac\-tible edge.
\end{theorem}

For $4$-connected graphs, the following theorem \cite[Theorem 45]{Kriesell2002}
can be deduced from Theorem \ref{T5} (sketch of proof see there).

\begin{theorem}
  \notinblk{\cite[Theorem 45]{Kriesell2002}}
  \label{Tred4}
  Every $4$-connected graph non-isomorphic to $K_5$ or $K_{2,2,2}$ can be reduced to a smaller
  $4$-connected graph by contracting one or two edges.
\end{theorem}

As it has been mentioned in the introduction, the option of contracting a constantly bounded number of edges in one step
might yield a similar reduction theorem for $5$-connected graphs, but not for $k$-connected graphs where $k>5$.

Theorems \ref{Tred2}, \ref{Tred3}, and \ref{Tred4} also provide the minor-minimal $k$-connected graphs
for $k \leq 4$. It follows from \cite[(1.4)]{RobertsonSeymour1990}
that the class of minor-minimal $5$-connected graphs is finite (as there exist planar $5$-connected graphs),
but  there is an `exact' conjecture by Fijav\v{z} \cite{Fijavz2001}:

\begin{conjecture}
  \notinblk{\cite{Fijavz2001}}
  \label{C2}
  Every $5$-connected graph contains a minor isomorphic to
  one of the graphs $K_6$, $K_{2,2,2,1}$, $C_5 * \overline{K_3}$,
  $I$, $\tilde{I}$, or $G_0$.%
  \footnote{Here $K_6$ is the complete graph on six vertices, the
  Tur\'an--graph $K_{2,2,2,1}$
  is obtained from a complete graph on seven vertices by deleting
  three independent edges, $C_5*\overline{K_3}$ is obtained from
  a cycle $C_5$ by adding $3$ new vertices and making them adjacent
  to all vertices of the $C_5$,
  $I$ denotes the icosahedron,
  $\tilde{I}$ is the graph obtained from $I$
  by replacing the edges of a cycle $abcdea$ induced by
  the neighborhood of some vertex with the edges of a cycle $abceda$,
  and $G_0$ is the graph obtained from the icosahedron by
  deleting a vertex $w$, replacing the edge $ab$ of a cycle $abcdea$ induced
  by the neighborhood of $w$ with the two edges $ac$ and $ad$,
  and, finally, identifying $b$ and $e$.}
\end{conjecture}

Fijav\v{z} proved this for graphs embeddable on the projective plane
($K_{2,2,2,1}$, $C_5 * \overline{K_3}$ are not projective planar) \cite{Fijavz2001,Fijavz2001b}. See also Table \ref{TB2}.
Apparently, to prove that, for any $k>5$, there are only finitely many minor-minimal $k$-connected graphs,
we need the full statement of Wagner's Conjecture (proved in \cite{RobertsonSeymour2004}).

\stepcounter{footnote} \footnotetext{By a Theorem of Dirac \cite{Dirac1960},
every graph of minimum degree at least $3$ contains a subdivision of $K_4$.}

\begin{table}
\centering
\begin{minipage}{\textwidth}
\begin{tabular}{r|c|l}
  $k$ & Minor base & Reference \\ \hline
  $1$ & $\{K_2\}$ & obvious \\
  $2$ & $\{K_3\}$ & obvious \\
  $3$ & $\{K_4\}$ & well-known${}^{\thefootnote}$, Theorem \ref{Tred3}, Tutte \cite{Tutte1961} \\
  $4$ & $\{K_5,K_{2,2,2}\}$ & Theorem \ref{Tred4} \\
  $5$ & $\{K_6, K_{2,2,2,1}, C_5 * \overline{K_3}, I, \tilde{I}, G_0\}$ & conjectured by Fijav\v{z} \cite{Fijavz2001,Fijavz2001b}
\end{tabular}
\end{minipage}
\caption{\label{TB2}
Numbers $k$ for which the minor-minimal $k$-connected graphs are known or predicted.}
\end{table}

A different reduction method has been developed by Dawes in \cite{Dawes1983}. His starting point was
Dirac's Theorem that the class of minimally $2$-connected graphs is the class of graphs obtained from
the cycles by finite sequences of attaching paths of length at least $2$ to suitable pairs of vertices%
\footnote{The resulting graph is always $2$-connected, and it is not too difficult to characterise the {\em suitable}
pairs for which the result is {\em minimally} $2$-connected.}.
The corresponding reduction theorem thus states that every minimally $2$-connected graph
distinct from a cycle can be reduced to a smaller such graph by deleting the interior vertices of some path of length at least
$2$ whose interior vertices have degree $2$ (not every such path will do the trick).%
\footnote{As an easy consequence, we get the well-known theorem on ear-decompositions, that every $2$-connected graph
can be obtained from a cycle by subsequently attaching paths of length at least $1$ to suitable pairs of vertices:
If $G$ is any $2$-connected graph non-isomorphic to a cycle then its simple subdivision is minimally $2$-connected and
Dirac's Theorem gives us a path which corresponds to a path $P$ in $G$ such that deleting all its elements
except its endvertices produce a smaller $2$-connected graph.}
In \cite{Dawes1983}, Dawes suggested the following operations for each $k \geq 1$ to construct a larger
graph from a given one:
{\em Operation $A_k$} is to choose $s \geq 1$ distinct edges and $k-2s \geq 0$ distinct
vertices, delete each of the chosen edges, add a new vertex $x$, and add a new edge from $x$ to each
end vertex of each chosen edge and from $x$ to each
chosen vertex. Operation $A_k$ is a special Henneberg construction.
{\em Operation $B_k$} is to choose $s \geq 2$ distinct edges and $k-s-1 \geq 0$ distinct vertices,
subdivide each of the chosen vertices, and add an edge from every subdivision vertex
to every other subdivision vertex and to every chosen vertex.
{\em Operation $C_k$} is to choose $k$ vertices, add a new vertex $x$, and add an edge from $x$ to every chosen
vertex. Here is Dawes's construction theorem for minimally $3$-connected graphs \cite[Theorem 6]{Dawes1983} \cite{Dawes1986}:

\begin{theorem}
  \notinblk{\cite[Theorem 6]{Dawes1983} \cite{Dawes1986}}
  \label{TDawes}
  The class of minimally $3$-connected graphs is the class of graphs obtained from $K_4$ by
  finite sequences of operation $A_3,B_3,C_3$ at suitable sets of objects.\footnote{Again, the graph resulting
  from applying $A_3$, $B_3$, or $C_3$ to {\em any} set of objects meeting the definition of the respective operation
  is $3$-connected --- but here it is much more difficult to characterise those sets which yield a {\em minimally}
  $3$-connected graph. See \cite{Dawes1983,Dawes1986} for the details.}
\end{theorem}

Theorem \ref{TDawes} implies the following reduction theorem for minimally $3$-connected graphs \cite{Dawes1983,Dawes1986}:

\begin{theorem}
  \notinblk{\cite{Dawes1983,Dawes1986}}
  Every minimally $3$-connected graph non-isomorphic to $K_4$ can be reduced to a smaller minimally $3$-connected
  graph by homotopically deleting edges or deleting vertices.
\end{theorem}

Theorem \ref{TDawes} `literally' holds for minimally $1$- and $2$-connected graphs; however,
note that the preconditions of $A_1$, $B_1$, and $B_2$ cannot be fulfilled, so that we get back the well-known fact that
the minimally $1$-connected graphs are the non-trivial trees, and that the class minimally $2$-connected graphs
is the class obtained from $K_3$ by finite sequences of subdivisions ($A_2$) and attaching paths of length $2$
(which follows from Dirac's theorem). Dawes conjectured in \cite[p. 287]{Dawes1983} that these facts
plus Theorem \ref{TDawes} generalise to generator theorems for minimally $k$-connected graphs for each $k \geq 4$
but also noticed that the `compatibility condition' characterizing the suitable sets of objects would become
much more difficult. --- However, let us disprove his conjecture:

\begin{theorem}
  For each $k \geq 4$, there are infinitely many minimally $k$-connected graphs which do
  not arise from a smaller minimally $k$-connected graph by operation $A_k$, $B_k$, or $C_k$.
\end{theorem}

{\bf Proof.}
First note that $A_k$ generates a $k$-connected graph only if all chosen edges are independent
and not incident with the chosen vertices (otherwise the new vertex will have less than $k$ neighbors).
$B_k$ generates a $k$-connected graph only if the chosen edges are not incident with the chosen vertices
(otherwise some subdivision vertex will have less than $k$ neighbors).
For all three operations, the degree at the endvertices of the chosen edges
does not change after application, and the degree of each chosen vertex
increases by at least one. Moreover, $B_k$ generates graphs with triangles unless
$s+k-s-1 \leq 2$, so unless $k \leq 3$.

For $k \geq 4$ it follows that if a
{\em $k$-regular} $k$-connected graph
arises from some (minimally) $k$-connected graph $G$
by applying $A_k$, $B_k$, or $C_k$
then it must arise by applications of either $A_k$ to $G$ and $t=\frac{k}{2}$ independent edges where $k$ is even,
or $B_k$ with $s=k-1$ distinct edges (no vertices are chosen).

Consequently, for odd $k \geq 5$, every $k$-regular $k$-connected graph
which arises from some (minimally) $k$-connected graph by $A_k,B_k,C_k$ has to contain a triangle ---
but there are infinitely many triangle free $k$-regular $k$-connected graphs.

For even $k \geq 4$, let the graph $H$ arise from a multicycle
of length at least $5$ with edge multiplicity $\frac{k}{2}$ by
subdividing each edge once. Every edge of $H$ is incident with
precisely $k$ others, two edges of $H$ can not be separated
by removing less than $k$ edges, and $H$ is edge-transitive. Its line graph
$G=L(H)$ is, consequently, $k$-regular, $k$-connected, and vertex-transitive.
We can not obtain it from applying $A_k$ (with $s=k/2$) to some minimally $k$-connected graph,
since the neighborhood of the new vertex as in $A_k$
had to contain at least $\frac{k}{2}$ pairwise disjoint pairs
of non-adjacent vertices --- but this situation does not occur at
any vertex of $G$.
Up to isomorphism, there is a unique graph $G^-$ from which $G$
arises by applying $B_k$ (with $s=k-1$),
but $G^-$ has a separating vertex set of order $\frac{k}{2}+1<k$.
($G^-$ can be obtained from $G$ by 
contracting any complete subgraph $K_k$ to a single vertex.)

Therefore, for even $k \geq 4$, too, there are infinitely many minimally $k$-connected graphs which do
not arise from a smaller one by $A_k,B_k,C_k$.
\hspace*{\fill}$\Box$

Yet another approach to constructively characterise the $k$-connected graphs
is Slater's concept of `splitting and soldering' \cite{Slater1974,Slater1976,Slater1978}.
Let $x$ be a vertex of a graph $G$ and $A,B$ sets such
that $A \cup B = N_G(x)$. To {\em split $x$ into $(a,b)$ according to $(A,B)$} means to delete $x$ from $G$,
add two new vertices $a,b$ and an edge from $a$ to $b$, and add an edge from $a$ to every $y \in A$
and from $b$ to every $z \in B$. If $|A|,|B| \geq k$ and $|A \cap B|=0$ then we say that the
new graph arises from $G$ by {\em $k$-vertex-splitting}, and if $|A|,|B| \geq k$ and $|A \cap B|=1$
then we say that the new graph arises from $G$ by {\em $k$-edge-splitting}%
\footnote{The `edge' to which this notion refers is the one connecting $x$ to the vertex in $A \cap B$ in $G$.}.
In this terminology, for example, the class of $3$-connected (simple) graphs is the
class of graphs obtained from $K_4$ by finite sequences of edge addition, $3$-vertex-splitting, and
$3$-edge-splitting --- as a constructive version of the Wheel Theorem (Theorem \ref{Twheel}).

Let $G$ be a graph and $K$ be a complete subgraph of order $k$ of $G$.
To {\em solder $x$ on $K$} means to add a new vertex $x$ and an edges from $x$ to every $y \in K$ to $G$,
and to delete a certain set $F$ of edges from $E(K)$. If, in the new graph $H$, $d_H(x) \geq k$ for all $x \in V(K)$
and $(V(K),F)=\overline{H[V(K)]}$ does not contain a $4$-cycle then we say that $H$ arises from $G$
by {\em $k$-soldering}. The reason for excluding $C_4$ here is to guarantee that $H$ is $k$-connected
if $G$ is, which is otherwise not true in general \cite{Slater1978}. The drawback of this operation is that $G$ might contain
a minor which is not a minor of $H$.

One of the main results in \cite{Slater1974} is that the class of $3$-connected (simple) graphs
is the class of graphs obtained from $K_4$ by finite sequences of line addition and $3$-soldering.
This does not similarly extend to $4$-connected graphs:
$4$-connected line graphs of cubic graphs (see above) are examples that cannot be obtained
from any smaller $4$-connected graph by edge-addition, $4$-vertex-splitting, $4$-edge-splitting, or $4$-soldering.
Slater's solution to overcome this situation was to generalise the vertex splitting as follows.
Let $x$ be a vertex of a graph $G$ and $A_1,\dots,A_r$ sets such that $A_1 \cup \dots \cup A_r=N_G(x)$.
To {\em split $x$ into $(a_1,\dots,a_r)$ according to $(A_1,\dots,A_r)$} means to delete $x$ from $G$,
add $r$ new vertices $a_1,\dots,a_r)$ and all edges connecting any two of them, and add an edge
from $a_j$ to every $y \in A_j$, $j \in \{1,\dots,r\}$. If $|A_j| \geq k$ for every $j \in \{1,\dots,r\}$ and
the $A_j$ are pairwise disjoint then we say that the new graph arises from $G$ by {\em $r$-fold $k$-vertex-splitting}.
The following is the second main result from \cite{Slater1974}.

\begin{theorem}
  \notinblk{\cite{Slater1974}}
  The class of $4$-connected graphs is the class of graphs obtained from $K_5$ by finite
  sequences of edge addition, $4$-soldering, $4$-ver\-tex-split\-ting, $4$-edge-splitting,
  and $3$-fold $4$-vertex-splitting. 
\end{theorem}

A more recent reduction theorem for $4$-connected graphs, due to Saito, uses a special Henneberg reduction;
these reductions have been very successfully applied in the context of edge-connectivity and aboricity questions.
Let $x$ be a vertex of some graph $G$ and let $\sigma$ be a partition of $E_G(x)$ into classes of cardinality $1$ or $2$.
The graph $G \stackrel{\sigma}{-} x$ is the graph obtained from $G-x$ by adding a new edge%
\footnote{That might be a loop.} between the endvertices distinct from $x$ of each pair of edges
$e \not= f$ such that $\{e,f\} \in \sigma$, and we say that it arises from $G$
by a {\em Henneberg reduction of degree $|\sigma|$ at $x$}.
The nice feature when performing a Henneberg reduction of degree $\ell$ at a vertex of degree $2\ell$
is that the degree function of the result is equal to the degree function of $G$ on $V(G)-\{x\}$;
in particular, when applied to a $k$-connected graph, the resulting graph will have minimum degree at least $k$.
The drawback is, as for soldering, that $G \stackrel{\sigma}{-} x$ might contain a minor which is not a minor of  $G$,
and conversely. Here is Saito's Theorem (from \cite[Theorem 2]{Saito2003}):

\begin{theorem}
  \notinblk{(See \cite[Theorem 2]{Saito2003}.)}
  For every vertex $x$ of a $4$-connected graph $G \not\cong K_5$,
  there exists a $4$-contractible edge at distance at most $1$ from $x$, or
  there exists a vertex $y$ of degree $4$ at distance $1$ from $x$
  and a partition $\sigma$ of $E_G(y)$ into two sets of cardinality $2$ such that
  $G \stackrel{\sigma}{-} y$ is $4$-connected.  
\end{theorem}

Another attempt to reduce a given $3$-connected graph $G$ to a smaller one is to {\em homotopically delete a vertex $x$},
i. e. to delete $x$ and then repeatedly suppress vertices of degree at most $2$ as long as it is possible.
This might do severe damage to $G$: We could kill the entire graph, for example if $G$ is a wheel and $x$ is its center.
The resulting graph is, however, well-defined \cite{Kriesell2008}, and denoted by $G -- x$. Let us mention the main
result from \cite{Kriesell2008}:

\begin{theorem}
  \notinblk{\cite{Kriesell2008}}
  Every $3$-connected graph non-isomorphic to $K_{3,3}$, $K_2 \times K_3$, or a wheel
  can be reduced to a smaller $3$-connected graph by homotopically deleting a vertex.
\end{theorem}

Although homotopic vertex deletion is by no means a `bounded' operation, 
it still has a constructive counterpart in terms of series parallel extensions \cite{Kriesell2008}.
As the homotopic deletion of vertices or edges can be considered as a (not constantly bounded)
sequence of edge deletions, vertex deletions, or edge contractions, it still would keep us in
any class closed under taking minors.

Here is an example of a contraction/deletion result for the class of triangle free $3$-connected graphs from \cite{Kriesell2007c}:

\begin{theorem}
  \notinblk{\cite{Kriesell2007c}}
  Every $3$-connected triangle free graph non-isomorphic to $K_{3,3}$ or the skeleton of a $3$-dimensional cube
  can be reduced to a smaller $3$-connected triangle free graph by  at most six edge deletions, vertex deletions,
  or edge contractions.
\end{theorem}

Let us finally mention a reduction theorem for bipartite graphs. Apart from trivial situations,
we would leave the class when {\em contracting} an edge in a bipartite graph.
Instead of identifying a pair of non-adjacent vertices, I have suggested to identify
a pair of (non-adjacent) vertices of the same colour class, which led indeed to a reduction theorem
\cite[Theorem 4]{Kriesell1999}:

\begin{theorem}
  \notinblk{\cite[Theorem 4]{Kriesell1999}}
  \label{Tbip}
  Every $k$-connected bipartite graph on more than $2k^2-2k+2$ vertices
  and, for even $k$, non-isomorphic to a graph $C_\ell[\overline{K_{k/2}}]$ for some even $\ell \geq 4$,
  can be reduced to a smaller $k$-connected bipartite graph by identifying two distinct non-adjacent vertices.
\end{theorem}

The bound to the number of vertices is sharp for infinitely many $k$, as it is shown by the point-line-incidence graphs
of projective geometries \cite{Kriesell1999}. It is easy to see that the `irreducible' bipartite graphs in Theorem \ref{Tbip}
for $k=2$ are the cycles of even length, and in \cite[Theorem 5]{Kriesell1999}
the nine irreducible graphs for $k=3$ have been determined.
It might well be that there are similar results for the class of $k$-connected $3$-colourable graphs
or even generalizations to $r$-colourable graphs where $r>3$.

\section{Subgraph Deletion}

\label{Shcks}

All the previously considered local reduction problems can be reformulated as questions for the existence
of a subgraph $X$ in some $k$-connected graph $G$ such that $\kappa(G-V(X)) > k-|V(X)|$
(alternatively: $\kappa(G-E(X))>k-|E(X)|$), where, in addition, $X$ respects some size- or connectivity conditions,
or, more generally, $X$ belongs to a certain class ${\cal H}$ of graphs.
For example, when asking for a $k$-contractible edge in a graph $G$ we ask for a subgraph $X \cong K_2$ of $G$
such that $\kappa(G-V(X))>k-2$, and ${\cal H}$ would be the class consisting of $K_2$.

As we have seen, such an $X$ need not to exist. Is there one which keeps the connectivity high?
Let us formulate this more precisely. Given a class ${\cal H}$ of graphs,
does there exist a function $f(k)$ such that every $f(k)$-connected graph $G$ admits a subgraph $X$ from ${\cal H}$
such that $G-V(X)$ (alternatively: $G-E(X)$) is $k$-connected?

For every non-empty finite graph class ${\cal H}$, such an $f$
obviously exists (take $f(k):=k+\max\{|V(X)|:\,X \in {\cal H}\}$).
Another obvious case in which $f$ exists is that ${\cal H}$ contains a graph $X$ which is contained
in every sufficiently highly connected graph%
\footnote{Any forest $X$ is such a candidate, see, for example \cite{Bollobas1978}.}.

Apart from this, there are also classes ${\cal H}$ for which the answer is not obvious. 
Thomassen proved \cite{Thomassen1981} that  every $(k+3)$-connected graph $G$ has an induced cycle $C$ such that
$G-V(C)$ is $k$-connected (see next section); that is, $f(k):=k+3$ will give a `yes' to our question if
${\cal H}$ is the class of all cycles. On the other hand, not every infinite class ${\cal H}$ admits
such an $f$, simply because there might be highly connected graphs which do not contain
objects from ${\cal H}$ as subgraphs (let ${\cal H}$ be, for example, the odd cycles).
Therefore, a first step towards our question could be the more fundamental question of
{\em characterizing the graph classes ${\cal H}$ such that every
sufficiently highly connected graph contains a member from ${\cal H}$ as a subgraph}.
However, for these classes our problem is solved, as it has been pointed out by K\"uhn and Osthus \cite[p. 30]{KuehnOsthus2003}:

\begin{theorem}
  \notinblk{\cite[p. 30]{KuehnOsthus2003}}
  \label{TKH}
  Suppose that ${\cal H}$ is a class of graphs such that there is an $\ell$ such that every $\ell$-connected
  graph contains a member from ${\cal H}$ as a subgraph. Then, for every $k$, there exists an $f_{\cal H}(k)$ such
  that every $f_{\cal H}(k)$-connected graph contains a member $H$ of ${\cal H}$ such that
  $G-V(H)$ (alternatively: $G-E(X)$) is $k$-connected.
\end{theorem}

This follows from Theorem 1 in \cite{KuehnOsthus2003} that
the vertex set of every $2^{11} 3 k^2$-connected graph admits a partition into two sets $A,B$ such
that $G[A],G[B]$ are $k$-connected and every vertex in $A$ has at least $k$ neighbors in $B$;
any copy of $H \in {\cal H}$ in $G[A]$ will do (cf. \cite{KuehnOsthus2003}).
  
Therefore, we propose to exclude those graphs without subgraphs from ${\cal H}$:

\begin{problem}
  \label{P1}
  Determine (the) classes of graphs ${\cal H}$ for which there is a high-connectivity-keeping-${\cal H}$-theorem, i.e.:
  For every $k$ there exists an $f(k)$
  such that every $f(k)$-connected graph $G$ which contains a subgraph from ${\cal H}$
  admits a subgraph $X$ from ${\cal H}$ such that $G-V(X)$ (alternatively: $G-E(X)$) is $k$-connected.
\end{problem}

Table \ref{TB1} summarises the graph classes
${\cal H}$ for which a high-connectivity-keeping-${\cal H}$-theorem 
in the sense of Problem \ref{P1} exists or is conjectured.

\stepcounter{footnote} \footnotetext{A {\em theta graph} is any subdivision of the graph $K_4^-$,
the complete graph on $4$ vertices minus a single edge.}

\begin{small}
\begin{table}
\begin{center}
\begin{tabular}{r|c|l}
  ${\cal H}$ & $f(k)$ & Reference \\ \hline
  cycles & k+3 & Thomassen \cite{Thomassen1981} \\
  even cycles & k+4 & Fujita and Kawarabayashi \cite{FujitaKawarabayashiPreprint} \\
  odd cycles & ? & conj. Thomassen \cite{Thomassen2001} \\
  theta graphs${}^{\thefootnote}$ & k+4 & Fujita and Kawarabayashi \cite{FujitaKawarabayashiPreprint} \\
  $t$-connected graphs & $4k+4t-13$ & $t \geq 3$ Hajnal \cite{Hajnal1983} \\
  & $k+t+1$ ? & conj. Thomassen \cite[p. 167]{Thomassen1983} \\
  as in Theorem \ref{TKH} & $2^{11} 3 \max\{k^2,\ell^2\}$ & K\"uhn and Osthus \cite{KuehnOsthus2003}
\end{tabular}
\caption{\label{TB1}
Graph classes ${\cal H}$ for which a non-obvious high-connectivity-keeping ${\cal H}$-theorem as in Problem \ref{P1} exists or is conjectured.}
\end{center}
\end{table}
\end{small}

There is a number of conjectures and results of the same flavour.
The first one leads us back to contractibility. 
By Theorem \ref{Tred3}, every $3$-connected graph non-isomorphic to $K_4$
contains a connected subgraph $H$ on two vertices such that $G-V(H)$ is $2$-connected or, equivalently,
that the graph obtained from $G-V(H)$ by adding a new vertex $h$ and making it adjacent to all neighbors
of $V(H)$ in $G$ is $3$-connected. McCuaig and Ota conjectured the following generalization of this \cite{McCuaigOta1994}:

\begin{conjecture}
  \notinblk{\cite{McCuaigOta1994}}
  \label{Cmco}
  For every $\ell \geq 1$ there exists an $f(\ell)$ such that every $3$-con\-nec\-ted graph on at least $f(\ell)$ vertices
  contains a connected subgraph $H$ on $\ell$ vertices such that $G-V(H)$ is $2$-connected.
\end{conjecture}

Conjecture \ref{Cmco} is true for $k \leq 4$, where the optimal values are $f(2)=5$ \cite{Tutte1961},
$f(3)=9$ \cite{McCuaigOta1994}, and $f(4)=8$ \cite{Kriesell2000b}.
These values are sharp, and they are not monotone in $k$.
Moreover, every {\em cubic} $3$-connected graph on at least $13$ vertices has a contractible subgraph
on $5$ vertices, and it is maybe interesting to see that the `local character' of the proof yields a generalization
to $3$-connected graphs of average degree at most $3+1/132$ \cite{Kriesell2007b}.
The conjecture is wide open in general. For example, it is not known whether there exists a $k$
such that its restriction to $k$-connected graphs is true. However, if we do not insist that $H$ as
in the statement is {\em connected}, then there is a positive result \cite{Kriesell2001b}:

\begin{theorem}
  \notinblk{\cite{Kriesell2001b}}
  \label{Tmco}
  For every $\ell \geq 1$ there exists an $f(\ell)$ such that every $3$-connected graph on at least $f(\ell)$ vertices
  contains a subgraph $X$ on $\ell$ vertices such that $G-V(X)$ is $2$-connected.  
\end{theorem}

If we delete edges instead of vertices then there is the following positive result from \cite{Kriesell2007b}
(based on a Theorem of Lemos and Oxley \cite{LemosOxley2003}):
For every $\ell \geq 1$ there exists an $f(\ell)$ such that every $4$-connected graph on at least $f(\ell)$ vertices
contains a path or a star $X$ on $\ell$ vertices such that $G-E(X)$ is $2$-connected.  

Mader generalised the statement of Theorem \ref{Tmco} to higher connectivity, as follows \cite{Mader2004}:

\begin{theorem}
  \notinblk{\cite{Mader2004}}
  For every $k \geq 4$ and every $\ell \geq 2$ there exists an $f_k(\ell)$ such that every $k$-connected graph on
  at least $f_k(\ell)$ vertices contains a subgraph $X$ on $\ell$ vertices such that $G-V(X)$ is $(k-2)$-connected.
\end{theorem}

It is not possible to replace $k-2$ by $k-1$ here \cite{Mader2004}.

The possibly most prominent conjecture along these lines is Lov\'asz's Conjecture on high-connectivity keeping paths
(cited according to \cite[p. 267]{Thomassen1983b}):

\begin{conjecture}
  \notinblk{(cf. \cite[p. 267]{Thomassen1983b})}
  \label{Chckp}
  For every $k$ there exists an $f(k)$ such that for any two vertices $a,b$ of any
  $f(k)$-connected graph $G$, there exists an induced $a,b$-path $P$ in $G$ such
  that $G-V(P)$ is $k$-connected.
\end{conjecture}

This conjecture has been verified for $k \leq 2$, where $f(1)=3$ and $f(2)=5$ are best possible
(see \cite{ChenGouldYu2003}, \cite{Kriesell2001c}, and also \cite{KawarabayashiLeeYu2005}).
Some years ago, I have suggested to first prove a version where we delete $E(P)$ instead of $V(P)$
and/or do not insist that $P$ is induced. There is recent progress on this, due to Kawarabayashi, Lee, Reed, and Wollan
\cite{KawarabayashiLeeReedWollan2008}:

\begin{theorem}
  \notinblk{\cite{KawarabayashiLeeReedWollan2008}}
  \label{Thckp}
  For every $k$ there exists an $f(k)$ such that for any two vertices $a,b$ of any
  $f(k)$-connected graph $G$, there exists an $a,b$-path $P$ in $G$ such
  that $G-E(P)$ is $k$-connected.
\end{theorem}

It is easy to see that an affirmative answer to Conjecture \ref{Chckp} would imply Theorem \ref{Thckp},
even with the additional constraint that $P$ is induced. It is not clear if the `non-induced' version of
Conjecture \ref{Chckp} would imply Theorem \ref{Thckp}. It might be the case that Theorem \ref{Thckp}
generalises as follows:

\begin{conjecture}
  \label{Chckst}
  For every $k$ there exists an $f(k)$ such that every $f(k)$-con\-nec\-ted graph
  contains a spanning tree $T$ such that $G-E(T)$ is $k$-connected.
\end{conjecture}

The edge-connectivity of this is obviously true, because every $(2k+2)$-edge-connected graph
admits $k+1$ edge-disjoint spanning trees, by a well-known corollary of the base packing theorem
by Tutte and Nash-Williams; removing one of them yields a supergraph of the union of $k$ edge-disjoint
spanning trees and hence a $k$-connected graph. In particular, Conjecture \ref{Chckst} is true for $k=1$,
where $f(1)=4$ is best possible. Jord\'an proved that every $6k$-connected graph $G$ has
$k$ edge-disjoint $2$-connected subgraphs, which shows that $f(2) \leq 12$ \cite{Jordan2005}.
The edge sets of these subgraphs are actually bases of the {\em $2$-dimensional rigidity matroid} of $G$.
It might be that there is a way of proving conjecture \ref{Chckst} for larger $k$ by using properties
of the higher dimensional rigidity matroids of the graph in question, but these are objects which are
far from being well understood \cite{Jordan}.

\section{Partitions under Connectivity Constraints}

\label{Spart}

As an alternative to the notion of high-connectivity-keeping subgraphs,
most of the previous results can be considered as partition statements:
Given a sufficiently highly connected graph $G$ with certain extra properties,
we look for a partition $\{A,B\}$ of $V(G)$ such that $G[A]$ meets a `small' or `bounded'
size condition and, in many cases, a mild connectivity condition,
whereas we want $G[B]$ highly connected. In some cases, an additional condition
to the location of $A$ might be incorporated. For example, the presence of a contractible
edge in a graph of connectivity $k$ [incident with some vertex $x$]
is equivalent to the presence of such a partition, where $G[A]$ has order two and is
connected [and contains $x$], whereas $G[B]$ needs to be
$(k-1)$-connected. 

In this section we look at a problem where the conditions to
the partition sets are more balanced. Let me first mention Gy\H{o}ry's classic
characterization of the $k$-connected graphs \cite{Gyoeri1978}.

\begin{theorem}
  \notinblk{\cite{Gyoeri1978}}
  A graph $G$ is $k$-connected if and only if
  for any $k$ distinct vertices $a_1,\dots,a_k$,
  and any $k$ positive integers $b_1,\dots,b_k$ with $\sum_{i=1}^k b_i = |V(G)|$,
  there exist disjoint sets $C_1,\dots,C_k \subseteq V(G)$ such
  that for all $i \in \{1,\dots,k\}$, $G[C_i]$ is a connected graph on $b_i$ vertices containing $a_i$.
\end{theorem}

There are versions for digraphs and also for edge-connectivity instead of connectivity (see \cite{Gyoeri1978}).
It is maybe surprising that this theorem has not been employed so far to the type
of problems we mentioned here; but, on the other hand, the same applies,
basically,  to Menger's theorem, as the vast majority of the arguments run
exclusively in terms of ${\cal T}(G)$. The reason might be that the non-trivial
part is the necessity of the partition conditon, i. e. we could apply the theorem
to the graph $G$ under consideration but possibly not to certify a certain
connectivity of some substructure.

The following question is due to Thomassen \cite{Thomassen1983} (cf. Table \ref{TB1}):

\begin{conjecture}
  \notinblk{\cite{Thomassen1983}}
  \label{CCT}
  For every $(s+t+1)$-connected graph $G$, there exists a partition $\{A,B\}$ of $V(G)$
  such that $G[A]$ is $s$-connected and $G[B]$ is $t$-connected.
\end{conjecture}

The qualitative part of this question has been settled: The conclusion holds
for $(4s+4t+1)$-connected graphs $G$ \cite{Thomassen1983}. The proof runs in three steps.
\begin{enumerate}
  \item[(i)]
    Since $G$ has minimum degree at least $4s+4t+1$, there exists
    a partition $\{A'',B''\}$ of $V(G)$ such that $G[A'']$ has minimum degree at least $4s$
    and $G[B'']$ has minimum degree at least $4t$, by a famous theorem of Stiebitz \cite{Stiebitz1996}.
    In particular, $G[A''],G[B'']$ have average degree at least $4s,4t$ respectively.
  \item[(ii)]
    By a Theorem of Mader \cite{Mader1972b} (see also \cite{Diestel2005}),
    the average degree bound from (i) ensures that
    $G[A'']$ has an $s$-connected subgraph and
    $G[B'']$ has a $t$-connected subgraph.
    That is, there exist disjoint subsets $A' \subseteq A''$
    and $B' \subseteq B''$ such that $G[A']$ is $s$-connected and $G[B']$ is $t$-connected.
  \item[(iii)]
    Since $G$ is $(s+t-1)$-connected, we can extend $A',B'$ to a partition of $\{A,B\}$ of $V(G)$
    as desired, that is, $A \supseteq A'$, $B \supseteq B'$, and $G[A]$ is $s$-connected
    and $G[B]$ is $t$-connected; this is due to a beautiful argument of Thomassen \cite{Thomassen1983}.
\end{enumerate}

(i) and (iii) show that the conclusion of Conjecture \ref{CCT} holds also under the weaker
assumption that $G$ is $(s+t-1)$-connected and has minimum degree $4s+4t+1$.
(iii) shows that for proving the conclusion of Conjecture \ref{CCT} it suffices to
find two disjoint subgraphs, which are $s$- and $t$-connected, respectively.
By careful consideration of the original bounds from Mader \cite{Mader1972b}, which are slightly better than $4k$,
Hajnal improved the bound to $4s+4t-13$ for all $s,t \geq 3$ \cite{Hajnal1983}.

There is also a version of Conjecture \ref{CCT} where we partition $E(G)$ instead of $V(G)$,
which has been posed by Mader in \cite{Mader1979} as an open problem:

\begin{problem}
  \notinblk{\cite{Mader1979}}
  \label{CWM}
  Given $s,t \geq 2$, does every $(s+t)$-connected graph admit a partition $\{A,B\}$ of $E(G)$
  such that the graph formed by the edges of $A$ is $s$-connected and the graph formed by those of $B$ is $t$-connected?
\end{problem}
Let me sketch a proof for the qualitative part of it:
Suppose that $G$ is a $(2s+4t)$-connected graph,
or just any $s$-connected graph of average degree at least $2s+4t$. Then $G$
contains a minimally $s$-connected spanning subgraph $H$. As we have seen
in Section \ref{Smin}, the average degree of $H$ is less than $2s$,
so that the average degree of $G-E(H)$ is at least $4t$. 
By Mader's Theorem, $G-E(H)$ contains a $t$-connected subgraph $T$.
Since $S:=G-E(T)$ is a supergraph of $H$, it follows that the edge sets of $S$ and $T$
partition $E(G)$ in the desired way. 

The first non-trivial case of Conjecture \ref{CCT} has  been settled by Thomassen \cite{Thomassen1981}:

\begin{theorem}
  \notinblk{\cite{Thomassen1981}}
  \label{TCT}
  Every $(k+3)$-connected graph $G$ contains an induced cycle $C$
  such that $G-V(C)$ is $k$-connected.
\end{theorem}

The proof is by induction on $V(G)$ of the stronger statement that there is an induced cycle $C$
such that every vertex not in $V(C)$ has at most $3$ neighbors in $V(C)$ (such that
$G-V(C)$ is $k$-connected). Obviously, every triangle would serve as such a $C$,
so that we may assume that $G$ is triangle free. By Theorem \ref{Ttriangle}, $G$ contains a contractible edge
$e$. From a cycle in $G/e$ with the desired properties it is the easy to obtain one in $G$.
--- Along these lines, let me also point to Mader's result that every $(k+2)$-connected graph contains
a cycle $C$ such that $G-E(C)$ is $k$-connected \cite[p. 190]{Mader1974}. In fact, $C$ can be taken as an induced cycle here.

It seems to be extremely difficult to prove local versions of this statement.
For example, the statement that `there is a function $g$ such that every for every edge $e$ of
every $g(k)$-connected graph $G$ there exists an induced cycle $C$ containing $e$
such that $G-V(C)$ is $k$-connected' is equivalent to Conjecture \ref{Chckp}, as has been
observed by Thomassen. Even the following, weaker problem, is open:

\begin{conjecture}
  \label{CMK}
  For every $k$, there exists an $h(k)$ such that every for every vertex $v$ of
  every $h(k)$-connected graph $G$ there exists an induced cycle $C$ containing $v$
  such that $G-V(C)$ is $k$-connected.
\end{conjecture}

It might be that partition problems as in Conjecture \ref{CCT} or \ref{CMK} are easier to solve on graphs with high girth.%
\footnote{The {\em girth} $g(G)$ of a graph $G$ is the length of a shortest cycle and $+\infty$ if $G$ is a forest.}
This is supported by the following facts conjectured by Thomassen and proved by Egawa \cite{Egawa1987} \cite{Egawa1998}:

\begin{theorem}
  \notinblk{\cite{Egawa1987} \cite{Egawa1998}}
  Let $k \geq 2$.
  \begin{enumerate}
    \item
      Every $(k+2)$-connected triangle free graph $G$ contains an induced cycle $C$
      such that $G-V(C)$ is $k$-connected.
    \item
      Every $(k+1)$-connected graph $G$ of girth at least $5$ contains an induced cycle $C$
      such that $G-V(C)$ is $k$-connected.
  \end{enumerate}
\end{theorem}

Let us look at the case that $s=t=k$ in Conjecture \ref{CCT}.
Hajnal's Theorem implies that the vertex set of every $(8k-13)$-connected graph $G$ has a partition
into two sets $A,B$ such that $G[A],G[B]$ are $k$-connected. This improves by almost a factor of $2$
if we restrict the statement to graphs of girth larger than $k$:
By \cite[Korollar 2]{Mader1972b}, every graph with $\delta(G) \geq 2k-2$ and girth $g(G)>k$ contains
a $k$-connected subgraph;
therefore, following the (i)-(ii)-(iii)-argument right after Conjecture \ref{CCT},
the vertex set of every $(4k-3)$-connected graph with $g(G)>k$ has a partition into
two sets $A,B$ such that $G[A],G[B]$ are $k$-connected.
In fact, for $k=3$ it is possible to use \cite[Satz 1]{Mader1972b}, with $n=3$ and $\nu=-1$, in part (ii) of the argument,
which then yields that the vertex set
of every $9$-connected graph without an induced subgraph $K_4^-$\hspace*{2mm}%
\footnote{$K_4^-$ is obtained from $K_4$ by deleting a single edge.}
has a partition into two sets $A,B$ such that $G[A],G[B]$ are $3$-connected.
Also for the (i)-part of the argument, there are improvements for graphs of high girth:
Whereas the general bound resulting from Stiebitz's Theorem is $s+t+1$,
the vertex set of every {\em triangle-free} graph of minimum degree at least $s+t$
can be partitioned into $A,B$ with $\delta(G[A]) \geq s$ and $\delta(G[B]) \geq t$ \cite{Kaneko1996},
and the bound $s+t$ can be improved once more to $s+t-1$ for graphs of girth at least $5$ \cite{Diwan2000}.
(As a consequence, the vertex set of every $8$-connected triangle-free graph 
can be partitioned into $A,B$ such that $G[A],G[B]$ are $3$-connected.)

Also for {\em large} graphs, the bound in Mader's Theorem in (ii) of the argument right after Theorem \ref{CCT} improves:
It follows from \cite[Satz 4]{Mader1972b} that every sufficiently large graph of average degree at least $(2+\sqrt{2}) \cdot k$ contains
a $k$-connected subgraph. As we have no control on the size of the partition classes when applying Stiebitz' result in (i) of the argument, 
this statement does not improve the bound of $4(s+t)+1$ for large graph immediately.
However, if it would be true that (*) for $n,s,t$ there exists a $f(n,s,t)$ such that
every graph $G$ of minimum degree $s+t+1$ with $|G| \geq f(n,s,t)$
admits a partition $\{A'',B''\}$ of $V(G)$ such that $G[A'']$ has average degree
at least $s$, $G[B'']$ has average degree at least $t$, and $|A''|,B''| \geq n$, then we get the following:

\begin{theorem}
  For $s,t,n$ there exists $f(s,t,n)$ such that
  for every $(2+\sqrt{2}) \cdot (s+t)+1$-connected graph $G$ with $|G| \geq f(s,t,n)$, there exists a partition $\{A,B\}$ of $V(G)$
  such that $G[A]$ is $s$-connected and $G[B]$ is $t$-connected.
\end{theorem}

A (probabilistic) proof of (*) has been announced by Carmesin and et al. \cite{Carmesinetal2010}.

\section{Line Graphs}
\label{Slgr}

Most of the problems mentioned in Section \ref{Shcks} and \ref{Spart}
provide an (often literal) analogous version in terms of edge-connectivity
(instead of connectivity). The answers are often affirmative, mostly due to
the presence of enough edge-disjoint spanning trees. (See, for example,
the paragraph right after Conjecture \ref{Chckst}). Moreover, it is straightforward
to translate these results on graphs  into the language of their line graphs: For example,
the fact that for any two vertices $a,b$ of an $f(k)$-edge-connected graph
$G$ there exists an $a,b$-path $P$ such that $G-E(P)$ is $k$-edge-connected
immediately implies that for any two {\em edges} $e,f$ of an $(f(k)+2)$-edge-connected
graph $G$ there is a path with terminal edges $e,f$ such that $G-E(P)$ is $k$-connected; consequently,
in the line graph $L(G)$, between any two {\em vertices} there exists an {\em induced} path $P$
such that $L(G)-V(P)$ is $k$-connected. One is tempted to say that this proves Lov\'asz's Conjecture,
Conjecture \ref{Chckp}, for line graphs --- but this is not true, because a high-edge-connectivity of some graph
is sufficient but not necessary for high connectivity of its line graph\footnote{Think of pendant edges.}.
However, it is easy to see that if $L(G)$ is $k$-connected then the vertices of degree at least $k$ in $G$
are $k$-edge-connected in $G$ \cite{Kriesell2001c}. This observation has been combined with a powerful
theorem on removable paths in graphs with a given edge-connectivity function by Okamura \cite{Okamura1984},
in order to finally prove that Lov\'asz's Conjecture is true for line graphs \cite{Kriesell2001c}.

Let us give another example of how to employ this method. We prove Conjecture \ref{Chckst} for line graphs.

\begin{theorem}
  Every $(12k+11)$-connected line graph has a spanning tree $T$ such that $G-E(T)$ is $k$-connected.
\end{theorem}

{\bf Proof.}
Let $G$ be any graph such that $L(G)$ is $(6k+6)$-connected
and $L(G)$ has minimum degree at least $12k+11$.
By \cite[Lemma 1]{Kriesell2001c}, the set $A$ of vertices of degree at least $6k+6$ in $G$ is $(6k+6)$-connected,
and since $L(G)$ has minimum degree at least $12k+11$, $G-A$ is edgeless. 
By  \cite[Theorem 3.1]{FrankKiralyKriesell2003},
$G$ admits $2k+2$ edge-disjoint trees such that each of them covers $A$.
Therefore, there exists $k+1$ edge-disjoint connected subgraphs $G$ $T_1,\dots,T_k$ such
that each of them covers $A$, every vertex of $A$ has degree at least $2$ in every $T_j$,
and $E(T_1),\dots,E(T_{k+1})$ partitions $E(G)$.

For each vertex $x \in V(G)$, let $K_x$ be the clique on $E_G(x)$ in $L(G)$,
and for each $x \in V(G)$, let $K_x^j$ be the subclique induced by $E_{T_j}(x)$ in $K_x$.

For each $x \in A$, take $m_x$ such that $|K_x^m|=\max\{|K_x^{j}|:\,j \in \{1,\dots,k+1\}\}$.
Then $K_x^m$ is a clique on at least $4$ vertices and, thus, has a non-separating spanning path $M_x^{m_x}$.
For each $j \not= m_x$, there exists a matching $M_x^j$ in $K_x$ such that each edge of $M_x^j$
connects a vertex from $V(K_x^j)$ to a vertex from $V(K_x^{m_x})$ and each vertex of $V(K_x^j)$ is connected
this way. The graph $H_x^1$ formed by $\bigcup_{j=1}^{k+1}E(M_x^j) \cup E(K_x^1)$ is, therefore, a connected
spanning subgraph of $K_x$ such that, for each $j \not= 1$, the graphs $H_x^j:=K_x^j - E(H_x^1)$ are connected
and for every vertex $e \in V(K_x)$ and every $j$ such that $e \not\in V(H_x^j)$
there exists an edge in $E(K_x)-E(H_x^1)$ connecting $e$ to some vertex in $V(H_x^j)$.

Let $H^j:=\bigcup_{x \in V(G)-A} K_x^j \cup \bigcup_{x \in A} H_x^j$.

We claim that for $e \not= f \in L(G)$,
there exists an $e,f$-path in $H^1$ and $k$ openly disjoint $e,f$-paths in $L(G)-E(H^1)$ (*),
which proves the theorem.
There exists $x,y \in A$ such that $e \in V(K_x)$ and $f \in V(K_x)$, $e \in V(K_y)$.
If $x=y$ then there exists an $e,f$-path in $H^1$ as $H_x^1$ is a spanning connected subgraph of $K_x$.
Let $j_e,j_f \in \{1,\dots,k+1\}$ such that $e \in V(K_x^{j_e}),f \in V(K_x^{j_f})$. Take $s_{j_e}:=e,r_{j_f}:=f$,
for each $j \not= j_e$ take a vertex $s_j \in N_{K_x-E(H^1)}(e) \cap V(K_x^j)$, and
for each $j \not= j_f$ take a vertex $r_j \in N_{K_x-E(H^1)}(e) \cap V(K_x^j)$.
For $j \not=1$, there exists an $s_j,r_j$-path $L_j$ in $K_x^j-E(H_1)$, and as the
$L_j$ are disjoint by definition, the paths $eL_jf$ are openly disjoint $e,f$-paths, which proves (*).
If $x \not= y$ then let $P_j$ be any $x,y$-path in $T_j$, and let $Q_j$ be the path induced by $E(P_j)$ in $L(G)$.
By construction, $E(Q_1) \subseteq E(H^1)$ and $E(Q_j) \cap E(H^1) =\emptyset$ for all $j \not= 1$.
Since the $Q_j$ are edge disjoint, the $P_i$ are disjoint. Each $P_j$ connects a vertex $e_j$ from $V(K_x^j)$
to a vertex $f_j$ from $V(K_y^j)$. By construction, there exists an $e,e_j$-path $R_j$ such that $R_j-e$ is in $H_x^j$
and an $f_j,f$-path $S_j$ such that $S_j-f$ is in $H_y^j$. It follows that $S_1 Q_1 R_1$ is an $e,f$-path in $H^1$
and that $S_j Q_j R_j$ for $k \not= 1$ are $k$ openly disjoint paths in $L(G)-E(H^1)$.
\hspace*{\fill}$\Box$

Another example is that also Conjecture \ref{CCT} is true for line graphs, simply because a line graph
of minimum degree $s+t+1$ contains $K_{s+1}$ or $K_{t+1}$ as a subgraph.
Let me finally mention the main result from \cite{Kriesell1998} that there is no $3$-con-critically $k$-connected line graph
(which implies Slater's conjecture from \cite{MaurerSlater1977} that there is no $3$-critically $k$-connected line graph).

\begin{sloppy}
\bibliographystyle{plain}
\bibliography{dc_bibliography}
\end{sloppy}

\bigskip

\centerline{*}

{\bf Author's address.}

Matthias Kriesell \\
IMADA $\cdot$ SDU \\
Campusvej 55 \\
DK-5230 Odense M

Denmark

\end{document}